\documentclass[12pt,draftcls,onecolumn]{IEEEtran}
\usepackage{graphicx,epsfig,amssymb,cite}

\usepackage{graphicx,epsfig}

\date{}
\usepackage{mathrsfs,color}

\begin{document}
\title{ New Results on Stability of Singular Stochastic Markov Jump Systems with State-Dependent Noise }
\author{{\normalsize Yong Zhao ${}^{1,2}$} , Weihai Zhang ${}^{1,*} $\thanks{%
${}^*$ Corresponding author.   {\protect\small  E-mail address:
 w\_hzhang@163.com}} \\
%EndAName
${}^1$  {\small College of Electrical Engineering and Automation, }\\
\small {Shandong University of Science and Technology,
Qingdao 266590,  China}\\
${}^2$ {\small  College of Mathematics and Systems Science}, \\
\small {Shandong University of Science and Technology,Qingdao
266590,  China}}

 \maketitle
{\em Abstract-} {This paper aims to develop the stability theory for
singular stochastic Markov jump systems with state-dependent noise,
including  both  continuous- and discrete-time cases. The sufficient
conditions for the existence and uniqueness
 of a solution to the system equation are provided. Some new and fundamental concepts such as
 non-impulsiveness and  mean square admissibility  are introduced, which are different
from those of  other existing works. By making use of
 the $\mathcal{H}$-representation technique and  the pseudo inverse $E^+$ of a singular matrix $E$,
sufficient conditions ensuring the  system   to be mean square
admissible are established in terms of strict  linear matrix
inequalities, which can be regarded as extensions of  the
corresponding results of deterministic singular systems  and normal
 stochastic systems.
  Practical examples are given to demonstrate the
effectiveness of the proposed approaches.}

{\em   Index Terms-} { Singular stochastic systems; stability;
$\mathcal{H}$-representation; matrix inequalities.}

\section{Introduction}

Over the past years, a great deal of attention  has  been paid to
singular systems because of its extensive applications to many
practical problems, such as  electric circuits \cite{s35}, economics
[2-3], aircraft modeling [4], network theory [5-6], robotics
\cite{s39} and so on. To date, many significant notions and results
based on state-space systems, for instance, stability analysis
[8-12], robust control [9-13], $H_{\infty}$ control [12, 14],
$H_{\infty}$ filtering [12, 15], $H_2/H_{\infty}$ control [16], have
been successfully extended to singular systems.

On the other hand, in recent years, there have been increasing
interest in   the  singular hybrid  system driven by Markov jump
parameters as
\begin{equation}\label{eq:eq1}
E{\dot x}(t)=A(r(t))x(t), x_0\in R^{n}.
\end{equation}
 This is because system (1)
 can be used to model systems with structure change
resulting of some inner discrete event
  such as component failures or
 repairs, abrupt environmental disturbances [17-18].
 More efforts have been done on stability analysis  related to system (1), and a variety  of
 interesting results have been reported, we refer the reader to [17-22]
 and the references therein.

Nevertheless, it is noted that the system is often perturbed by
various types of environment noises in many branches of science and
industry [23-27]. If we take the environment noises into account,
system (\ref{eq:eq1}) becomes a singular stochastic differential
equation with Markov jump parameters, which can be represented as:
\begin{equation}\label{eq:eq2}
Edx(t)=A(r_t)x(t)dt+C(r_t)x(t)dw(t),
\end{equation}
where $w(t)$ is one-dimensional,
 standard Wiener process. Such systems are much more advanced and realistic,
 which possess all the characteristics of singular systems, It\^o equation, as well as hybrid systems.
  It should be pointed out that, so far,
   only few works were done on this class of systems as a result of  the influence  of the diffusion term.
    In fact, there are some difficulties in generalizing the results of deterministic
singular systems to singular It\^o systems. For instance, how do we
provide the condition
   for the existence and uniqueness of the solution to the system equation? how do we introduce
 the non-impulsiveness and admissibility of system (2)?
   For the cases of system (2) with or without Markov jump parameters, some articles either assumed that the solution of
the system equation existed  [28] or directly considered that the
regularity and non-impulsiveness of system (1) were the same as
system (2) [29-30], which completely ignored the crucial action of
the diffusion term $C(i)$. Also, as pointed out by [31], the
application of It\^o formula to system (2) without Markov jump
parameters was improper [28-29]. Recently, for system (2) in the
absence of
 Markov jump parameters, \cite{s25} laid a solid foundation for the further study of this
topic as It\^o formula was exactly applied to such systems; \cite{s26} further
investigated the observer-based controller design by using a
sequential design technique. In \cite{s27}, based on the
results of the singular Markov jump system (1),
 a sufficient condition for the mean square admissibility of system
(\ref{eq:eq2}) was given in the form of non-strict linear matrix
inequalities. This may cause a big trouble in checking the
conditions numerically.  Base on the above analysis, it can be found
that, up to now, the problem of stability analysis for
continuous-time singular stochastic  Markov jump systems with
state-dependent noise has not been fully investigated, even for the
discrete-time case, which motivates us to do this study.

In this paper, we are concerned with  the problem of stability
analysis for singular stochastic Markov jump systems with
state-dependent noise including the continuous- and discrete-time
cases.
 Conditions for the existence and uniqueness of the
solution to the system equation  are given. The new definitions for
non-impulsiveness and mean square admissibility of systems
considered are introduced, which are different from those in
existing works [28-35]. An equivalent conversion of singular
stochastic Markov jump systems with state-dependent noise into
standard and deterministic singular systems is obtained by using the
$\mathcal{H}$-representation technique, then the sufficient
condition for system considered to be mean square admissible is
presented in strict LMI, which is less conservative than that in
\cite{s27}. Furthermore, with the help of the the pseudo inverse
$E^+$ of the
 singular matrix $E$, new conditions for the mean square admissibility of the
continuous- and discrete-time singular stochastic system with state-dependent noise
 are formulated as strict LMI.
The results obtained  can be viewed  as extensions of the
corresponding results on stability of deterministic singular systems
and normal stochastic systems. Finally, illustrative examples are
given to show the effectiveness of the proposed results.

 Notations: $R^n$: space of all $n$- dimensional real vectors with usual $2$-norm $\|\cdot \|$;
 $R^{m\times n}$:  space of all $m\times n$ real matrices;
 $S_n$:  set of all $n\times n$ symmetric matrices;
 $S_n^{N}$: set of  $X=(X(1),\cdots,X(N))$ with $X(i)\in S_n, i=1,\cdots,N$;
 $A>0$(resp.$A<0$): $A$
 is a real symmetric positive definite (resp. negative definite) matrix;
  $A^T$: the transpose of $A$;  $vec(\cdot)$: the row stacking operator;
   $A\otimes B$: the Kronecker product of two matrices $A$ and $B$;
   $E^+$: the Moore-Penrose pseudo inverse  of the matrix $E$;
    $\mathbb{E}(\cdot)$: the expectation operator; $I_n$: the
    $n\times n$ identity matrix.
    $R(A)$: the rank of $A$.

%\hfill mds
%
%\hfill January 11, 2007
\section{Preliminaries}

Let us consider the following continuous-time singular stochastic
It\^o's system with Markov jump parameters:
\begin{equation}\label{eq:eq3}
\left\{
\begin{array}{lcl}
Edx(t)=A(r_t)x(t)dt+C(r_t)x(t)dw(t), \ \ \\
Ex(0)=x_{0}\in R^{n},\ \
\end{array}
\right.
\end{equation}
where $E$ is a constant $n\times n$-dimensional matrix with $rank(E)=r\leq n $; $x(t)\in R^{n}$ is the system state vector, $x_0\in R^{n} $ is the initial condition;
 $w(t)$ is one-dimensional,
 standard Wiener process that is defined on the filtered probability space $(\Omega, \mathcal F, {\mathcal F_t}, {\mathcal P})$ with a filtering $\{\mathcal F_t\}_{\{t\geq 0\}}$; $\{r_t, t\geq 0\}$ is a right continuous homogeneous Markov chain taking values in a finite state space $S=\{1,\cdots,N\}$ with transition probability matrix $\prod=\{\pi_{ij}\}_{N\times N}$ given by
 \begin{equation}\label{eq:eq4}
\pi_{ij}=p\{r_{t+h}=j|r_t=i\}=\left\{
\begin{array}{lcl}
\pi_{ij}h+o(h),   i\neq j, \\
1+\pi_{ii}h+o(h), i=j,
\end{array}
\right.
\end{equation}
 where $h>0$, $\lim_{h\rightarrow 0}{o(h)/h=0}$, and $\pi_{ij}\geq 0 (i\neq j)$ represents the transition rate from $i$ to $j$, which satisfies $\pi_{ii}=
 -\sum_{j\neq i}{\pi_{ij}}$.

For the sake of discussion, we introduce two operators $\phi:\ S_n^N\rightarrow R^{n^2N}$  and $\psi:\ S_n^N\rightarrow  R^{\frac{n(n+1)N}{2}}$ which are defined as follows:
\begin{equation}\label{eq:eq5}
\phi(X)=\left[\begin{array}{cc}vec(X(1))\\\cdots\\vec(X(N))\end{array}\right],
\psi(X)=\left[\begin{array}{cc}svec(X(1))\\\cdots\\svec(X(N))\end{array}\right],
\forall X=(X(1),\cdots,X(N))\in S_n^N,
\end{equation}
where $vec(X(i))=(x_{11}(i),\cdots,x_{1n}(i),\cdots,x_{n1}(i),\cdots,x_{nn}(i))^T\in R^{n^2}, svec(X(i))=(x_{11}(i),\\
\cdots,x_{1n}(i),x_{22}(i),\cdots,
x_{2n}(i),\cdots,x_{n-1,n-1}(i),x_{n-1,n}(i),x_{nn}(i))^T\in R^{\frac{n(n+1)}{2}}$,
$x_{ks}(i)(k,s=1,\\2,
\cdots,n)$ is the element of $X(i)$. Note that $svec(X(i))$  is derived by deleting the repeated
elements of $vec(X(i))$. By Definition 2.1 of \cite{s20}, the relation between $svec(X(i))$ and $vec(X(i))$ is given via $vec(X(i))=H_nsvec(X(i))$,
 where $H_n$ is a $n^2\times \frac{n(n+1)}{2}$-dimensional matrix independent of $X(i)$.\\
By Lemmas 2.1, 2.3 and 2.4 of \cite{s20}, the following result can
be easily obtained.

\noindent{\bf Lemma 1:} $\mathrm{(i)}$ For any $X\in S_n^N$, there is a $n^2N\times{\frac{n(n+1)}{2}N}$ matrix $H_n^N$
 independent of $X$ such that $\phi(X)=H_n^N\psi(X)$, where $H_n^N=diag(H_n,\cdots,H_n)$.
 Conversely, for any $Y\in R^{n^2N\times {\frac{n(n+1)}{2}N}}$, there is $X\in{S_n^N}$ such that $\phi(X)=H_n^NY$,
 $\mathrm{(ii)}$ $(H_n^N)^TH_n^N$ is nonsingular, i.e., $H_n^N$ has full column rank.

\noindent{\bf Lemma 2:} Let $Z_1,Z_2\in S_n^N$, if $(H_n^N)^T\phi(Z_1)=(H_n^N)^T\phi(Z_2)$, then $\phi(Z_1)=\phi(Z_2)$.

{\bf Proof:} By Lemma 1-$\mathrm{(i)}$, $\phi(Z_1)=H_n^N\psi (Z_1)$, $\phi(Z_2)=H_n^N\psi (Z_2)$. From $(H_n^N)^T\phi(Z_1)=(H_n^N)^T\phi(Z_2)$,
 we have
 $$(H_n^N)^TH_n^N\psi (Z_1)=(H_n^N)^TH_n^N\psi (Z_2).$$
 By Lemma 1-$\mathrm{(ii)}$,  the above formula yields $\psi (Z_1)=\psi (Z_2)$, which is equivalent to $\phi(Z_1)=\phi(Z_2)$.

The next lemma plays a crucial role in this paper, which guarantees
the existence and uniqueness of the solution to system (3).

\noindent{\bf Lemma 3:} If there exist a pair of  nonsingular matrices
$M(i)\in R^{n\times n}$ and $ N(i)\in R^{n\times n}$ for every $i\in S$ such that one of the following conditions is satisfied,
then system  (\ref{eq:eq3}) has a unique solution.

{\small \begin{eqnarray}\label{eq:eq6}
(i)&&M(i)EN(i)=\left[\begin{array}{cc}I_{n_1}(i)&0\\0&N_{n_2}(i)\end{array}\right],
 M(i)A(i)N(i)=\left[\begin{array}{cc}\tilde A_{11}(i)&0\\0&I_{n_2}(i)\end{array}\right],\nonumber\\
 &&M(i)C(i)N(i)=\left[\begin{array}{cc}\tilde C_{11}(i)&\tilde C_{12}(i)\\0&0\end{array}\right],
 \end{eqnarray}}
 where $N_{n_2}(i)\in R^{n_2\times n_2}$ is a nilpotent, $N^h_{n_2}(i)=0$,
 $\tilde A_{11}(i),\tilde C_{11}(i)\in R^{n_1\times n_1},
 \tilde C_{12}(i)\in R^{n_1\times n_2},n_1+n_2=n$.
{\small\begin{eqnarray}\label{eq:eq7}
\hspace{-3mm}(ii)^{[34]} && M(i)EN(i)=\left[\begin{array}{cc}I_r(i)&0\\0&0\end{array}\right],  M(i)A(i)N(i)=\left[\begin{array}{cc}\bar A_{11}(i)&0\\0&I_{n-r}(i)\end{array}\right],\nonumber \\
 \hspace{-3mm}&&M(i)C(i)N(i)=\left[\begin{array}{cc}\bar C_{11}(i)&\bar C_{12}(i)\\0&\bar C_{22}(i)\end{array}\right],
 \end{eqnarray}}

 where $\bar A_{11}(i),\bar C_{11}\in R^{r\times r}$, $\bar C_{12}(i)\in R^{r\times (n-r)}$ , $\bar C_{22}(i)\in R^{(n-r)\times (n-r)}$.

{\bf Proof:} The proof of item (ii) can be found in Lemma 2.2 of
\cite{s27} . As to (i), let
$\xi(t)=N(i)^{-1}x(t)=[\xi_1(t)^T\quad\xi_2(t)^T]^T,\xi_1(t)\in
R^{n_1},\xi_2(t)\in R^{n_2}$, then under condition (i), system (3)
is equivalent to
\begin{equation}\label{eq:eq8}
d\xi_1(t)=\tilde A_{11}(i)\xi_1(t)dt+(\tilde C_{11}(i)\xi_1(t)+\tilde C_{12}(i)\xi_2(t))dw(t)
\end{equation}
and
\begin{equation}\label{eq:eq9}
N_{n_2}(i)d\xi_2(t)=\xi_2(t)dt.
\end{equation}
Take the Laplace  transform  on both sides of (\ref{eq:eq9}), we
have
\begin{equation}\label{eq:eq10}
(sN_{n_2}(i)-I)\xi_2(s)=N_{n_2}(i)\xi_2(0).
\end{equation}
From (\ref{eq:eq10}), we obtain
\begin{equation}\label{eq:eq11}
\xi_2(s)=(sN_{n_2}(i)-I)^{-1}N_{n_2}(i)\xi_2(0).
\end{equation}
The inverse Laplace transform of $\xi_2(s)$ yields
\begin{equation}\label{eq:eq12}
\xi_2(t)=-\sum_{i=1}^{h-1}\delta^{i-1}(t)N_{n_2}(i)\xi_2(0),
\end{equation}
where $L[\delta^{i}(t)]=s^i$.

After substituting (\ref{eq:eq12}) into (\ref{eq:eq8}),
(\ref{eq:eq8}) becomes an ordinary stochastic differential equation.
By Theorem 3.8 of \cite{s17}, the solution of (\ref{eq:eq8}) exists
and is unique, so does (3).

{\bf Remark 1:} A new condition for the existence and uniqueness of
the solution to system (3) is presented in Lemma 3.
  When $N_{n_2}(i)=0$, condition (i) is exactly the same as Remark 2.1-(I) of \cite{s27}. When $C(i)=0$, condition (i)
reduces to the corresponding result for the singular Markov jump
system (see \cite{s10},\cite{s12}-\cite{s16}).

In what follows, we  introduce several concepts which are essential
in this paper.

{\bf Definition 1:} System (3) is said to be impulse free if Lemma 3-(i) with
$deg(det(sE-A(i)))=R(E)$ or Lemma 3-(ii) hold for  every $i\in S$.

 {\bf Definition 2:} System (3) is said to be asymptotically stable in the mean square
 if for any initial condition $x_0\in{R^n} $,  $\lim_{t \rightarrow \infty}\mathbb{E}\{\|x(t)\|^2|x_0,r_0\}=0$.

 {\bf Definition 3:} System (3) is said to be mean square admissible if
 it has a unique solution and is impulse-free and asymptotically stable  in the mean square.

{\bf Remark 2:} As it is observed that  Definition 1 for the
non-impulsiveness of singular Markov jump systems with
state-dependent noise includes the diffusion term $C(i)$, which is
firstly introduced and different from  all the existing works
(\cite{s21}-\cite{s28}). Also, it is  obvious that the
non-impulsiveness of
 system (3) implies that system (3) has a solution, which coincides  with
 deterministic singular systems; see, e.g.,  \cite{s2} and \cite{s7}.
In the cases of $C(i)=0$ and $S=\{1\}$, Definition 1 and Definition
3 degenerate  to the corresponding definitions of deterministic
singular systems (\cite{s1}-\cite{s11}).

Parallel to the results of deterministic singular systems on the
non-impulsiveness behavior, the following proposition is obtained.

{\bf Proposition 1} System (3) is said to be impulse free if
one of the following conditions holds.\\
(i) $deg(det(sE-A(i)))=R(E)$ in (\ref{eq:eq6}).\\
(ii)$N_{n_2}(i)=0$ in (\ref{eq:eq6}).\\
(iii)$ \bar A_{22}(i)$ in the following (13) is invertible  and $R(E,C(i))=R(E)$,
\begin{eqnarray}
&&M(i)EN(i)=\left[\begin{array}{cc}I_r(i)&0\\0&0\end{array}\right],
 M(i)A(i)N(i)=\left[\begin{array}{cc} \bar A_{11}(i)&\bar A_{12}(i)\\\bar A_{21}(i)&\bar A_{22}(i)\end{array}\right],\nonumber\\
&&M(i)C(i)N(i)=\left[\begin{array}{cc}\bar C_{11}(i)&\bar C_{12}(i)\\0&0\end{array}\right],
\end{eqnarray}
 where $M(i)$ and $N(i)$ are invertible matrices.

{\bf Proof:} We will prove that conditions (i)-(iii) are equivalent.\\
 $(i)\Rightarrow (ii)$. If (\ref{eq:eq6}) holds, then
$det(sE-A(i))=det(sI_{n_1}(i)-\tilde A_{11}(i))det(sN_{n_2}(i)-I_{n_2}(i))$. But
$deg(det(sE-A(i))=R(E)=r$ implies $N_{n_2}(i)=0$ and $n_1=r$. \\
$(ii)\Rightarrow (i)$. It is easy to see that  $N_{n_2}(i)=0$ in (\ref{eq:eq6}) implies $deg(det(sE-A(i)))=R(E)$.\\
$(iii)\Rightarrow (i)$. Since $R(E)=r$ and $R(E,C(i))=R(E)$, it is always possible to choose two nonsingular
matrices $M(i)$ and $N(i)$ such that (13) holds. If $\bar A_{22}(i)$ is invertible,
there are nonsingular matrices $\bar M(i)$ and $\bar N(i)$ such that
{\small\begin{eqnarray}
 &&\bar M(i)E\bar N(i)=\left[\begin{array}{cc}I_r(i)&0\\0&0\end{array}\right],
 \bar M(i)A(i)\bar N(i)=\left[\begin{array}{cc}\bar A_{11}(i)-\bar A_{12}(i)\bar A_{22}^{-1}(i)\bar
 A_{21}(i)&0\\0&I_{n-r}(i)\end{array}\right],\nonumber\\
  &&\bar M(i)C(i)\bar N(i)=\left[\begin{array}{cc}\bar C_{11}(i)-\bar C_{12}(i)\bar A_{22}^{-1}(i)\bar A_{21}(i)&\bar C_{12}(i)\\0&0\end{array}\right].
\end{eqnarray}}
Obviously, $deg(det(sE-A(i)))=R(E)$ in (\ref{eq:eq6}).\\
$(i)\Rightarrow (iii)$. By $(i)\Leftrightarrow (ii)$, $N_{n_2}(i)=0$ in (\ref{eq:eq6}) implies
 $R(E,C(i))=R(E)$. Furthermore, there are nonsingular matrices $M(i)$ and $N(i)$ such that (13) holds. As a result,
 $deg(det(sE-A(i)))=deg(det(sI_r(i)-\bar A_{11}(i))det(-\bar A_{22}(i))$. By $deg(det(sE-A(i))=R(E)$,
we have $det(\bar A_{22}(i))\neq 0$, i.e., $\bar A_{22}(i)$ is invertible.\\
The following lemma will be used later.

 {\bf Lemma 4\  \cite{s29}: }For any three matrices $A,B$ and $C$ of
suitable dimensions, $vec(ABC)=(A\otimes C^T)vec(B)$.

 \section{Stability of continuous-time systems}
 In this section, we will focus on studying the stability of system (3) via two different approaches.
  One is to apply the $\mathcal H$-representation technique in \cite{s20} to improve the existing results of
  \cite{s27},
   and the other one is to establish  strict LMI conditions for the stability of system (3).

{\bf Theorem 1:} System (3) is mean square admissible if the following system
\begin{equation}\label{eq:eq13}
\mathcal{E}\psi(\dot X(t))=\mathcal{A}\psi(X(t))
\end{equation}
is admissible,  where (\ref{eq:eq13}) is an
$\frac{n(n+1)}{2}N$-dimensional deterministic singular system,
\begin{eqnarray*}
&&\mathcal{E}=(H_n^N)^Tdiag(E\otimes E,\cdots,E\otimes E)H_n^N, \\
&&\mathcal{A}=(H_n^N)^T[diag(A(1)\otimes E+E\otimes A(1)+C(1)\otimes C(1),\cdots,A(N)\otimes E\\
&&+E\otimes A(N)+C(N)\otimes C(N))+(\Pi\otimes I_{n^2})diag(E\otimes E,\cdots,E\otimes E)]H_n^N,\\
&&X_i(t)=\mathbb{E}\{x(t)x(t)^T{\chi_{\{r_t=i\}}}\},
X(t)=(X_1(t),\cdots,X_N(t)),
\end{eqnarray*}
and $x(t)$ is the trajectory of (3).

{\bf Proof:} For (3), let
$X_i(t)=\mathbb{E}\{x(t)x(t)^T{\chi_{\{r_t=i\}}}\}$, then use the
generalized It\^o's formula \cite{s18}, we derive
\begin{equation}\label{eq:eq14}
E\dot X_i(t)E^T=A(i)X_i(t)E^T+EX_i(t)A(i)^T+C(i)X_i(t)C(i)^T+\sum_{j=1}^{N}\pi_{ij}EX_j(t)E^T.
\end{equation}
Take $vec(\cdot)$ on both sides of (\ref{eq:eq14}), we obtain
{\small\begin{equation}\label{eq:eq15}
(E\otimes E)vec(\dot
X_i(t))=(A(i)\otimes E+E\otimes A(i)+C(i)\otimes C(i))vec(X_i(t))
+\sum_{j=1}^{N}\pi_{ij}(E\otimes E)vec(X_j(t)).
\end{equation}}
(\ref{eq:eq15}) can be equivalently written as
{\small\begin{eqnarray}\label{eq:eq16}
&diag(E\otimes E,\cdots,E\otimes E)\phi(\dot X(t))=[diag(A(1)\otimes E+E\otimes A(1)+C(1)\otimes C(1),\cdots,\nonumber\\
&A(N)\otimes E+E\otimes A(N)+C(N)\otimes C(N))+(\Pi\otimes I_{n^2})diag(E\otimes E,\cdots,E\otimes E)]\phi(X(t)).
\end{eqnarray}}
Premultiply  (\ref{eq:eq16}) by $(H_n^N)^T$  , then by Lemma 1-(i),
we have {\small\begin{eqnarray}\label{eq:eq17}
&&(H_n^N)^Tdiag(E\otimes E,\cdots,E\otimes E)H_n^N\psi(\dot X(t))=(H_n^N)^T[diag(A(1)\otimes E+E\otimes A(1)+C(1)\otimes C(1),\nonumber\\
&&\cdots,A(N)\otimes E+E\otimes A(N)+C(N)\otimes C(N))+(\Pi\otimes I_{n^2})diag(E\otimes E,\cdots,E\otimes E)]H_n^N\psi(X(t)).
\end{eqnarray}}
Let
\begin{eqnarray*}
&&\mathcal{E}=(H_n^N)^Tdiag(E\otimes E,\cdots,E\otimes E)H_n^N, \\
&&\mathcal{A}=(H_n^N)^T[diag(A(1)\otimes E+E\otimes A(1)+C(1)\otimes C(1),\cdots,A(N)\otimes E+E\otimes A(N)\\
&&\quad\quad+C(N)\otimes C(N))+(\Pi\otimes I_{n^2})diag(E\otimes
E,\cdots,E\otimes E)]H_n^N,
\end{eqnarray*}
(\ref{eq:eq17}) is converted into (\ref{eq:eq13}). Since
\begin{eqnarray*}
&&\lim_{t \rightarrow \infty}\psi(X(t))=0\Leftrightarrow\lim_{t  \rightarrow\infty}\phi(X(t))=0\\
&&\Leftrightarrow\lim_{t \rightarrow \infty}X(t)=0\Leftrightarrow\lim_{t \rightarrow \infty}\mathbb{E}\{\|x(t)\|^2|x_0,r_0\}=0,
\end{eqnarray*}
if (\ref{eq:eq13}) is admissible, then system (3) is mean square
admissible . This completes the proof.

{\bf Remark 3:} It is shown that (\ref{eq:eq16}) is a deterministic
but not a  standard  singular system, because $\phi(X(t))$ contains
repeated or redundant components. Generally speaking, for such a
nonstandard system, some classical criterions do not hold
\cite{s20}, this is why we introduce Lemmas 1-2.

{\bf Remark 4:} Compared  with the equation  (9) of \cite{s27},
singular stochastic systems (3) in this paper is equivalently
converted to a deterministic  $\frac{n(n+1)}{2}N$-dimensional
``standard''
 singular system (15). In this case, we can directly
make use of the existing results of deterministic singular systems
such as Theorem 2.2 of [12]  to provide the strict LMI condition for
the mean square admissibility of singular stochastic system (3).

{\bf Corollary 1:} System (1) is asymptotically mean square
admissible if there exist matrices $\mathcal P>0$ and $\mathcal Q$
such that the following inequality
\begin{equation}\label{eq aa}
(\mathcal P\mathcal E+\mathcal S\mathcal Q)^T\mathcal A+\mathcal
A^T(\mathcal P\mathcal E+\mathcal S\mathcal Q)<0
\end{equation}
 holds,  where $\mathcal S$ is any matrix with full column rank and satisfies $\mathcal E^T\mathcal S=0$.

{\bf Remark 5:} Theorem 1 tells us that we may make use of the
existing results of system (\ref{eq:eq13}) to discuss the stability
of system (3) provided that   both (3) and (\ref{eq:eq13})  have
solutions.  Theorem 1  is not convenient to be used in
  studying
stabilization,   $H_\infty$  control and other topics  of system
(3), which motivates us  to develop other methods.

{\bf Assumption 1:} For every $i\in S$, $R(E,C(i))=R(E)$ holds.

{\bf Lemma 5:}  Consider system (3) with initial state $(x_0,i)\in
R^n\times S$. If $V(x(t),t,r_t)=x^T(t)E^TP(r_t)x(t)$ with
$E^TP(r_t)=P^T(r_t)E$, then for any $t\geq 0$,
{\small\begin{equation} \mathbb E\{\int_0^t
\mathcal{L}V(x(s),s,r_s)ds|r_0=i\} =\mathbb
E\{x^T(t)E^TP(r_t)x(t)|r_0=i\}-x^T(0)E^TP(i)x(0),
\end{equation}}
where
{\small \begin{eqnarray}\label{eq:eq18}
&&\mathcal{L}V(x(t),t,r_t)=x^T(t)[A^T(r_t)P(r_t)+P^T(r_t)A(r_t)\nonumber\\
&&\quad\quad\quad+\sum_{j=1}^{N}\pi_{r_tj}E^TP(j)+C^T(r_t)(E^{+})^TE^TP(r_t)E^{+}C(r_t)]x(t).
\end{eqnarray}}
{\bf Proof}: By generalized It\^o's formula in \cite{s18} and
Theorem 2 of \cite{s25}, the above lemma  can be verified.

{\bf Theorem 2:} System (3) is mean square admissible if there exist
matrices $P(i)$ for each $i\in S$ such that
\begin{equation}\label{eq:eq19}
E^TP(i)=P^T(i)E\geq 0,
\end{equation}
\begin{equation}\label{eq:eq20}
A^T(i)P(i)+P^T(i)A(i)+\sum_{j=1}^{N}\pi_{ij}E^TP(j)+C^T(i)(E^{+})^TE^TP(i)E^{+}C(i)<0.
\end{equation}
Proof:  (i) First, we will show that system (3) has a solution and
is impulse free.

Under  Assumption 1, there exist a pair of nonsingular matrices $M,
N$ such that
\begin{eqnarray}\label{eq:eq21}
 &&MEN=\left[\begin{array}{cc}I_r&0\\0&0\end{array}\right],
 MA(i)N=\left[\begin{array}{cc}\bar A_{11}(i)&\bar A_{12}(i)\\\bar A_{21}(i)&\bar A_{22}(i)\end{array}\right],\nonumber\\
  &&MC(i)N=\left[\begin{array}{cc}\bar C_{11}(i)&\bar C_{12}(i)\\0&0\end{array}\right],
\end{eqnarray}
where the partitions have the appropriate dimensions. From
(\ref{eq:eq19}) and (\ref{eq:eq20}), we derive
\begin{equation}\label{eq:eq22}
A^T(i)P(i)+P^T(i)A(i)+\sum_{j=1}^{N}\pi_{ij}E^TP(j)<0.
\end{equation}
(\ref{eq:eq19}) and  (\ref{eq:eq22}) show that $(E,A(i))$ is
regular, impulse free and $\bar A_{22}(i)$ is invertible (see
Theorem 10.1 of \cite{s7}). Let
\begin{eqnarray}\label{eq:eq23}
 \bar M(i)=\left[\begin{array}{cc}I&-\bar A_{12}(i)\bar A_{22}^{-1}(i)\\0&\bar
 A_{22}^{-1}(i)\end{array}\right]M, \
\bar N(i)=N\left[\begin{array}{cc}I&0\\-\bar A_{22}^{-1}(i)\bar A_{21}(i)&I\end{array}\right],
\end{eqnarray}
we have
\begin{eqnarray}\label{eq:eq24}
 &&\bar M(i)E\bar
 N(i)=\left[\begin{array}{cc}I_r&0\\0&0\end{array}\right], \
 \bar M(i)A(i)\bar N(i)=\left[\begin{array}{cc}\bar A_{11}(i)-\bar A_{12}(i)\bar A_{22}^{-1}(i)\bar
 A_{21}(i)&0\\0&I\end{array}\right],\nonumber\\
  &&\bar M(i)C(i)\bar N(i)=\left[\begin{array}{cc}\bar C_{11}(i)-\bar C_{12}(i)\bar A_{22}^{-1}(i)\bar A_{21}(i)&\bar C_{12}(i)\\0&0\end{array}\right].
\end{eqnarray}
 It is shown that (\ref{eq:eq24}) satisfies Proposition 1-(ii), so (3) has a solution and is impulse-free.

(ii) Below, we will prove that system (3) is asymptotically stable in the mean square.

 Let
 \begin{equation}\label{eq:eq25}
\bar N(i)^{-1}x(t)=[\xi_1(t)^T\quad\xi_2(t)^T]^T,
\end{equation}
where $\xi_1(t)\in R^r, \xi_2(t)\in R^{n-r}$, then system (3) is
equivalent  to
\begin{eqnarray}\label{eq:eq26}
\left\{ \begin{array}{ll}
d\xi_1(t)=(\bar A_{11}(i)-\bar A_{12}(i)\bar A_{22}^{-1}(i)\bar A_{21}(i))\xi_1(t)dt\\
\quad \quad \quad+(\bar C_{11}(i)-\bar C_{12}(i)\bar A_{22}^{-1}(i)\bar A_{21}(i))\xi_1(t)dw(t) \\
\xi_2(t)=0.
\end{array}\right.
\end{eqnarray}
Let
\begin{eqnarray}\label{eq:eq27}
\bar M^{-T}(i)P(i)\bar N^{-1}(i)=\left[\begin{array}{cc}
\bar P_{11}(i)& \bar P_{12}(i)\\ \bar P_{21}(i)&\bar P_{22}(i)\end{array}\right],
\end{eqnarray}
(\ref{eq:eq19}) and (\ref{eq:eq22}) imply $\bar P_{11}(i)=\bar
P_{11}^T(i)>0, \bar P_{12}(i)=0$. The  Lyapunov function candidate
is selected as
\begin{eqnarray}\label{eq:eq28}
V(x(t),t,r_t)=x^T(t)E^TP(r_t)x(t)=\xi_1^T(t)\bar
P_{11}(r_t)\xi_1(t).
\end{eqnarray}
Take the expectation on both sides of (\ref{eq:eq28}),  we obtain
\begin{equation}\label{eq:eq29}
\lambda_{min}(\bar P_{11}(r_t))\mathbb E(\|\xi_1(t)\|^2)\leq \mathbb E(V(x(t),t,r_t))\leq \lambda_{max}(\bar P_{11}(r_t))\mathbb E(\|\xi_1(t)\|^2).
\end{equation}
From (\ref{eq:eq18}) and (\ref{eq:eq20}), there exists a positive
constant $\theta$ such that
\begin{equation}\label{eq:eq32}
\mathbb E(\mathcal LV(x(t),t,r_t))\leq -\theta\mathbb E(\|x(t)\|^2).
\end{equation}
Using (\ref{eq:eq25}), we have
\begin{equation}\label{eq:eq33}
\lambda_{min}(\bar N^T(r_t)\bar N(r_t))\mathbb E(\|\xi_1(t)\|^2)\leq
\mathbb E(\|x(t)\|^2)\leq \lambda_{max}(\bar N^T(r_t)\bar
N(r_t))\mathbb E(\|\xi_1(t)\|^2).
\end{equation}
Therefore,
\begin{equation}\label{eq:eq34}
\mathbb E(\mathcal LV(x(t),t,r_t))\leq -\theta \lambda_{min}(\bar
N^T(r_t)N(r_t))\mathbb E(\|\xi_1(t)\|^2).
\end{equation}
We define a new function $U(t)=e^{\rho t}V(x(t),t,r_t)$, then
\begin{equation}\label{eq:eq30}
dU(t)=\rho e^{\rho t}V(x(t),t,r_t)dt+e^{\rho t}dV(x(t),t,r_t).
\end{equation}
By the theory of stochastic differential equations, $\mathbb
E(\|\xi_1(t)\|^2)$ is  continuous and  bounded on any finite horizon
$[0,t]$, so do  $\mathbb E(V(x(t),t,r_t))$ and $\mathbb E(\mathcal
LV(x(t),t,r_t))$  due to (\ref{eq:eq29}) and (\ref{eq:eq34}). Hence,
we can proceed to  integrate  and take  the expectation on both
sides of (\ref{eq:eq30}), which  yields
\begin{equation}\label{eq:eq31}
\mathbb E(U(t))=\mathbb E(U(0))+\int_{0}^{t}\rho e^{\rho s}\mathbb E(V(x(s),s,r_s))ds+\int_{0}^{t}e^{\rho s}\mathbb E(\mathcal LV(x(s),s,r_s))ds.
\end{equation}

Substitute (\ref{eq:eq29}) and (\ref{eq:eq34}) into (\ref{eq:eq31}),
we derive
 \begin{eqnarray}\label{eq:eq35}
min_{r_t\in S}\{\lambda_{min}\hspace{-6mm}&&(\bar P_{11}(r_t))\}e^{\rho t}\mathbb E(\|\xi_1(t)\|^2)
\leq \int_{0}^{t}e^{\rho s}\{\rho max_{r_s\in S}\{\lambda_{max}\bar P_{11}(r_s)\}\}\mathbb E(\|\xi_1(s)\|^2)ds\nonumber\\
\hspace{-6mm}&&-\int_{0}^{t}e^{\rho s}\{\theta min_{r_s\in S}\{\lambda_{min}(N^T(r_s)N(r_s)\}\}\mathbb E(\|\xi_1(s)\|^2)ds+\mathbb E(U(0)).
\end{eqnarray}
 Select an appropriate $\rho$ which satisfies  $\rho \leq \frac {\theta min_{r_t\in S}\{ \lambda_{min}(N^T(r_t)N(r_t))\}}{max_{r_t\in S}\{\lambda_{max}\bar P_{11}(r_t)\}}$, we have
 \begin{eqnarray}\label{eq:eq36}
&&\mathbb E(\|\xi_1(t)\|^2)\leq \frac{\mathbb E(U(0))}{min_{r_t\in S}\{ \lambda_{min}\bar P_{11}(r_t)\}}e^{-\rho t}
=\frac{\xi^T_1(0)\bar P_{11}(r_0)\xi_1(0)}{max_{r_t\in S}\{ \lambda_{min}\bar P_{11}(r_t)\}}e^{-\rho t}\nonumber\\
&&\quad\quad\quad\quad\quad\leq \frac{max_{r_0\in
S}\{\lambda_{max}\bar P_{11}(r_0)\}}{min_{r_t\in S}\{
\lambda_{min}\bar P_{11}(r_t)\}}e^{-\rho t}\mathbb E(\|
\xi_1(0)\|^2).
\end{eqnarray}
Let $\gamma=\frac{max_{r_0\in S}\{\lambda_{max}\bar
P_{11}(r_0)\}}{min_{r_t\in S}\{\lambda_{min}\bar P_{11}(r_t)\}}$,
 (\ref{eq:eq36}) becomes
\begin{equation}\label{eq:eq37}
\mathbb E(\|\xi_1(t)\|^2)\leq \gamma e^{-\rho t}\mathbb E(\|
\xi_1(0)\|^2).
\end{equation}
Take  limit on (\ref{eq:eq37}), we obtain
\begin{equation}\label{eq:eq38}
\lim_{t \rightarrow \infty}\mathbb{E}\{\|\xi_1(t)\|^2\}=0.
\end{equation}
By Definition 2, ({\ref{eq:eq26}) is asymptotically stable in mean
square sense, so does (3). According to Definition 3, system (3) is
mean square admissible.

Note that (\ref{eq:eq19}) contains equality constraints, which are
usually not satisfied perfectly when solving them via LMI toolbox.
To this end, the following theorem provides a sufficient condition
for system (3) to be mean square admissible in terms of  strict
LMIs.

{\bf Theorem 3:} System (3) is mean square admissible if there exist
matrices   $P(i)>0$, $Q(i)$, for each $i\in S$
 such that
{\small \begin{equation}\label{eq:eq39}
A^T(i)(P(i)E+FQ(i))+(P(i)E+FQ(i))^TA(i)+\sum_{j=1}^{N}\pi_{ij}E^TP(j)E+C^T(i)(E^+)^TE^TP(i)EE^+C(i)<0,
\end{equation}}
where $F\in R^{n\times (n-r)}$ is any matrix with full column rank and satisfies $E^TF=0$.

{\bf Proof:} Let $\bar P(i)=P(i)E+FQ(i)$, it is easy to verify that
\begin{equation}\label{eq:eq40}
E^T\bar P(i)=\bar P^T(i)E=E^TP(i)E\geq 0,
\end{equation}
\begin{equation}\label{eq:eq41}
A^T(i)\bar P(i)+\bar P^T(i)A(i)+\sum_{j=1}^{N}\pi_{ij}E^T\bar  P(j)+C^T(i)(E^+)^TE^T\bar P(i)E^+C(i)<0.
\end{equation}
By Theorem 2, system (3) is mean square admissible.

{\bf Remark 6:}  When $E=I$, obviously, $F=0$. In this case,
Theorem 3 still coincides with the  stability theory of \cite{s19}
on state-space stochastic systems. When $S=\{1\}$, i.e., there is no Markov jump parameters in system (3),
Theorem 2 accords with Theorem 1 of \cite{s26}. When $C(i)=0$, Theorem 3
degenerates to  the criterion for the stability of singular Markov
jumping systems in \cite{s15}.

{\bf Remark 7:} It is noted that, under the assumption
$R(E,C(i))=R(E)$, Theorem 2 provides the sufficient conditions for
singular stochastic systems  to be mean square admissible. The
assumption condition is less conservative than all the previous
works (see, e.g. \cite{s25}-\cite{s28}), where Theorem 6.1 of
\cite{s28} was based  on the assumption that
$deg(det(sE-A))=R(E),det(sE-A)\neq 0,$ and $R(E,C(i))=R(E)$, while
Theorem 1 of \cite{s26} just omitted $deg(det(sE-A))=R(E)$. On the
other hand,
  since Theorem 2 includes
equality constraints, this may lead to computational  problem in
testing them. To overcome  this difficulty, the sufficient
conditions for the mean square admissibility of system (3) is
formulated in terms of strict LMI in Theorem 3. The strict LMI
conditions obtained are completely different from those in
\cite{s22}, \cite{s23},\cite{s26}-\cite{s28}, which  are  not only
easily to be tested by Matlab toolbox but also provide a useful way
to investigate other control problems such as state feedback,
$H_\infty$ control and $H_2/H_\infty$ control that have not been
fully discussed yet.

 \section{Stability of discrete-time systems}
Consider a class of discrete-time singular stochastic systems with
Markov jump parameters:
\begin{equation}\label{eq:eq42}
\left\{
\begin{array}{lcl}
Ex(k+1)=A(r_k)x(k)+C(r_k)x(k)w(k), \ \ \\
x(0)=x_{0}\in R^{n},\ \
\end{array}
\right.
\end{equation}
where $x(t)\in R^{n}$ is the system state, $x_0\in R^{n} $ is the
initial condition; $E$ is a constant $n\times n$-dimensional matrix
with $rank(E)=r\leq n $; $w(k)\in R$ is  a wide stationary,
second-order  process, $\mathbb{E}(w(k))=0$ and
$\mathbb{E}(w(k)w(s))=\delta_{ks}$ with $\delta_{ks} $ being a
Kronecker delta; the parameter $r_k$ represents a discrete-time
Markov chain taking values in a finite set $S=\{1,\cdots,N\}$ with
transition probabilities $Pr\{r_{k+1}=j|r_k=i\}=\lambda_{ij}$, and
the  transition probability matrix is given as
$\bigwedge=\{\lambda_{ij}\}_{N\times N}$, where $\lambda_{ij}\geq 0$
and satisfies $\sum_{j=1}^{N}\lambda_{ij}=1$ for any $i\in S$.

  To avoid  tedious repetition, results for discrete-time systems are only listed without any proof unless necessary.

{\bf Lemma 6:} If there exist a pair of  nonsingular matrices
$M(i)\in R^{n\times n}$ and $ N(i)\in R^{n\times n}$ for every $i\in
S$ such that one of the expressions (\ref{eq:eq6}), (\ref{eq:eq7})
holds, then (\ref{eq:eq42}) has a unique solution.

Parallel to the continuous-time case, the following definitions for
discrete-time stochastic singular system (\ref{eq:eq42}) are
presented.

 {\bf Definition 4:} System (\ref{eq:eq42}) is said to be causal if Lemma 3-(i) with
$deg(det(sE-A(i)))=R(E)$ or Lemma 3-(ii) holds for  every $i\in S$.

 {\bf Definition 5:} System (\ref{eq:eq42}) is said to be asymptotically stable in the mean square
 if it has a unique solution and for any initial condition $x_0\in{R^n} $, $\lim_{k \rightarrow \infty}\mathbb{E}\{\|x(k)\|^2|x_0,r_0\}=0$.

 {\bf Definition 6:} System (\ref{eq:eq42}) is said to be mean square admissible if
 it has a unique solution and is casual and asymptotically stable in the mean square.

Following the similar lines as those presented in Theorem 1, we can
obtain the following result readily.

{\bf Theorem 4:} System (\ref{eq:eq42}) is mean square admissible if
the following system
\begin{equation}\label{eq:eq43}
\mathcal{E}\psi(\dot X(k))=\tilde \mathcal{A}\psi(X(k))
\end{equation}
is admissible,  where (\ref{eq:eq43}) is an
$\frac{n(n+1)}{2}N$-dimensional deterministic singular system,
\begin{eqnarray*}
&&\mathcal{E}=(H_n^N)^Tdiag(E\otimes E,\cdots,E\otimes E)H_n^N, \\
&&\tilde \mathcal{A}=(H_n^N)^T(\Lambda^T\otimes I_{n^2})diag(A(1)\otimes A(1)+C(1)\otimes C(1),\cdots,
A(N)\otimes A(N)\\
&&+C(N)\otimes C(N))H_n^N,\\
&&X_i(k)=\mathbb{E}\{x(k)x(k)^T{\chi_{\{r_k=i\}}}\},
X(k)=(X_1(k),\cdots,X_N(t)).
\end{eqnarray*}

Now, We present new sufficient conditions for the admissibility of
system (\ref{eq:eq42}) in the form of LMIs.

{\bf Theorem 5:} System (\ref{eq:eq42}) is mean square admissible if
there exist symmetric matrices $P(i)=P^T(i)$ for each $i\in S$ such
that
\begin{equation}\label{eq:eq44}
E^TP(i)E\geq 0,
\end{equation}
\begin{equation}\label{eq:eq45}
A^T(i)(\sum_{j=1}^{N}\lambda_{ij}P(j))A(i)+C^T(i)(\sum_{j=1}^{N}\lambda_{ij}P(j))C(i)-E^TP(i)E<0.
\end{equation}

{\bf Proof:} (i) Our first aim is to prove that (\ref{eq:eq42}) has
a solution and is casual.

Based on  Assumption 1, there exist nonsingular matrices $M,N$ such
that
\begin{eqnarray}\label{eq:eq46}
 &&MEN=\left[\begin{array}{cc}I_r&0\\0&0\end{array}\right],
 MA(i)N=\left[\begin{array}{cc}\hat A_{11}(i)&\hat A_{12}(i)\\\hat A_{21}(i)&\hat A_{22}(i)\end{array}\right],\nonumber\\
  &&MC(i)N=\left[\begin{array}{cc}\hat C_{11}(i)&\hat C_{12}(i)\\0&0\end{array}\right].
\end{eqnarray}
Let
\begin{equation}\label{eq:eq47}
M^{-T}P(i)M^{-1}=\left[\begin{array}{cc}\hat P_{11}(i)&\hat P_{12}(i)\\\hat P^T_{12}(i)&\hat P_{22}(i)\end{array}\right],
\end{equation}
substitute (\ref{eq:eq46}) and (\ref{eq:eq47}) into (\ref{eq:eq44}),
(\ref{eq:eq45}) respectively, we derive
\begin{eqnarray}\label{eq:eq48}
\hspace{-15mm}&&E^TP(i)E=(N^{-T}\left[\begin{array}{cc}I_r&0\\0&0\end{array}\right]M^{-T})(M^T\left[\begin{array}{cc}\hat P_{11}(i)&\hat P_{12}(i)\\\hat P^T_{12}(i)&\hat P_{22}(i)\end{array}\right]M)(M^{-1}\left[\begin{array}{cc}I_r(i)&0\\0&0\end{array}\right]N^{-1})\nonumber\\
\hspace{-15mm}&&\quad \quad \quad \quad =N^{-T}\left[\begin{array}{cc}\hat P_{11}(i)&0\\0&0\end{array}\right]N^{-1}\geq 0,\end{eqnarray}
\begin{eqnarray}\label{eq:eq49}
\hspace{-15mm}&&A^T(i)(\sum_{j=1}^{N}\lambda_{ij}P(j))A(i)+C^T(i)(\sum_{j=1}^{N}\lambda_{ij}P(j))C(i)-E^TP(i)E\nonumber\\
\hspace{-15mm}&&=N^{-T}\left[\begin{array}{cc}W_1(i)&W_2(i)\\W^T_2(i)&W_3(i)\end{array}\right]N^{-1}<0,
\end{eqnarray}

where
{\small\begin{eqnarray*}
W_1(i)\hspace{-6mm}&&=\hat A^T_{11}(i)(\sum_{j=1}^{N}\lambda_{ij}\hat P_{11}(j))\hat A_{11}(i)+\hat C^T_{11}(i)(\sum_{j=1}^{N}\lambda_{ij}\hat P_{11}(j))\hat C_{11}(i)+\hat A^T_{21}(i)(\sum_{j=1}^{N}\lambda_{ij}\hat P^T_{12}(j))\hat A_{11}(i)\nonumber\\
\hspace{-6mm}&&+\hat A^T_{11}(i)(\sum_{j=1}^{N}\lambda_{ij}\hat P_{12}(j))\hat A_{21}(i)+\hat A^T_{21}(i)(\sum_{j=1}^{N}\lambda_{ij}\hat P_{22}(j))\hat A_{21}(i)-\hat P_{11}(i),\nonumber\\
W_2(i)\hspace{-6mm}&&=\hat A^T_{11}(i)(\sum_{j=1}^{N}\lambda_{ij}\hat P_{11}(j))\hat A_{12}(i)+\hat C^T_{11}(i)(\sum_{j=1}^{N}\lambda_{ij}\hat P_{11}(j))\hat C_{12}(i)+\hat A^T_{21}(i)(\sum_{j=1}^{N}\lambda_{ij}\hat P^T_{12}(j))\hat A_{12}(i)\nonumber\\
\hspace{-6mm}&&+\hat A^T_{11}(i)(\sum_{j=1}^{N}\lambda_{ij}\hat P_{12}(j))\hat A_{22}(i)+\hat A^T_{21}(i)(\sum_{j=1}^{N}\lambda_{ij}\hat P_{22}(j))\hat A_{22}(i),\nonumber\\
W_3(i)\hspace{-6mm}&&=\hat A^T_{12}(i)(\sum_{j=1}^{N}\lambda_{ij}\hat P_{11}(j))\hat A_{12}(i)+\hat C^T_{12}(i)(\sum_{j=1}^{N}\lambda_{ij}\hat P_{11}(j))\hat C_{12}(i)+\hat A^T_{22}(i)(\sum_{j=1}^{N}\lambda_{ij}\hat P^T_{12}(j))\hat A_{12}(i)\nonumber\\
\hspace{-6mm}&&+\hat A^T_{12}(i)(\sum_{j=1}^{N}\lambda_{ij}\hat P_{12}(j))\hat A_{22}(i)+\hat A^T_{22}(i)(\sum_{j=1}^{N}\lambda_{ij}\hat P_{22}(j))\hat A_{22}(i).\\
\end{eqnarray*}}
(\ref{eq:eq49}) implies $W_3(i)<0$. Furthermore, since $\hat
P_{11}(j)\geq 0$ for each $j\in S$, we have

\begin{equation}\label{eq:eq50}
\hat A^T_{12}(i)(\sum_{j=1}^{N}\lambda_{ij}\hat P_{11}(j))\hat A_{12}(i)+\hat C^T_{12}(i)(\sum_{j=1}^{N}\lambda_{ij}\hat P_{11}(j))\hat C_{12}(i)\geq 0.
\end{equation}
Therefore, in $W_3(i)$, it is easy to see that
{\small\begin{equation}\label{eq:eq51} \hat
A^T_{22}(i)(\sum_{j=1}^{N}\lambda_{ij}\hat P^T_{12}(j))\hat
A_{12}(i) +\hat A^T_{12}(i)(\sum_{j=1}^{N}\lambda_{ij}\hat
P_{12}(j))\hat A_{22}(i)+\hat
A^T_{22}(i)(\sum_{j=1}^{N}\lambda_{ij}\hat P_{22}(j))\hat
A_{22}(i)<0.
\end{equation}}
From ({\ref{eq:eq51}), it follows that $\hat A_{22}(i)$ is
invertible. By Definition 10.2 of \cite{s7}, $(E,A(i))$ is regular
and casual. Let
\begin{equation}\label{eq:eq52}
\hat N(i)=N\left[\begin{array}{cc}I_r(i)&0\\-\hat A^{-1}_{22}(i)\hat A_{21}(i)&\hat A^{-1}_{22}(i)\end{array}\right],
\end{equation}
we derive
\begin{eqnarray}\label{eq:eq53}
 && ME\hat N(i)=\left[\begin{array}{cc}I_r&0\\0&0\end{array}\right],
 MA\hat N(i)=\left[\begin{array}{cc}\hat A_{11}(i)-\hat A_{12}(i)\hat A_{22}^{-1}(i)\hat
 A_{21}(i)&\hat A_{12}(i)\hat A^{-1}_{22}(i)\\0&I_{n-r}(i)\end{array}\right],\nonumber\\
  &&MC(i)\hat N(i)=\left[\begin{array}{cc}\hat C_{11}(i)-\hat C_{12}(i)\hat A_{22}^{-1}(i)\hat A_{21}(i)&\hat C_{12}(i)\hat A^{-1}_{22}(i)\\0&0\end{array}\right].
\end{eqnarray}
By Proposition 1-(ii), system ({\ref{eq:eq42}) has a solution and is casual.\\
(ii) Now, we are in  a  position to prove that system
({\ref{eq:eq42}) is asymptotically stable in  mean square sense.
\\ Let
\begin{equation}\label{eq:eq54}
x(k)=\hat N(i)[\xi_1(k)^T\quad\xi_2(k)^T]^T,
\end{equation}
where $\xi_1(k)\in R^r, \xi_2(k)\in R^{n-r}$, then system ({\ref{eq:eq42}) is
equivalent  to
\begin{eqnarray}\label{eq:eq55}
\left\{ \begin{array}{ll}
\xi_1(k+1)=(\hat A_{11}(i)-\hat A_{12}(i)\hat A_{22}^{-1}(i)\hat A_{21}(i))\xi_1(k)\\
\quad \quad \quad \quad+(\hat C_{11}(i)-\hat C_{12}(i)\hat A_{22}^{-1}(i)\hat A_{21}(i))\xi_1(k)w(k) \\
\xi_2(k)=0.
\end{array}\right.
\end{eqnarray}
For simplicity, write
\begin{equation}\label{eq:eq56}
\tilde A_1(i)=\hat A_{11}(i)-\hat A_{12}(i)\hat A_{22}^{-1}(i)\hat A_{21}(i), \tilde C_1(i)=\hat C_{11}(i)-\hat C_{12}(i)\hat A_{22}^{-1}(i)\hat A_{21}(i).
\end{equation}
 Substitute (\ref{eq:eq47}) and (\ref{eq:eq53}) into (\ref{eq:eq45}), we obtain
\begin{eqnarray}\label{eq:eq57}
-\hat P_{11}(i)+\tilde A_1^T(i)(\sum_{j=1}^{N}\lambda_{ij}\hat P_{11}(j))\tilde A_1(i)
+\tilde C_1^T(i)(\sum_{j=1}^{N}\lambda_{ij}\hat P_{11}(j))\tilde C_1(i)<0.
\end{eqnarray}
By (\ref{eq:eq48}) and (\ref{eq:eq57}), it follows that $\hat
P_{11}(i)>0$ for every $i\in S$. Define
\begin{equation}\label{eq:eq58}
V(x(k),r_k)=x^T(k)E^TP(r_k)Ex(k)=\xi_1^T(k)\hat P_{11}(r_k)\xi_1(k)>0,
\end{equation}
it can be verified that
\begin{eqnarray}\label{eq:eq59}
\hspace{-10mm}&&\mathbb E\{V(\xi_1(k+1),r_{k+1}|\xi_1(k),r_k=i\}-V(\xi_1(k),r_k=i)\nonumber\\
\hspace{-10mm}&&=\xi^T_1(k)[-\hat P_{11}(i)+\tilde
A_1^T(i)(\sum_{j=1}^{N}\lambda_{ij}\hat P_{11}(j))\tilde A_1(i)
+\tilde C_1^T(i)(\sum_{j=1}^{N}\lambda_{ij}\hat P_{11}(j))\tilde
C_1(i)]\xi^T_1(k)\nonumber\\
&&<0.
\end{eqnarray}
Thus, there exists a positive constant $\delta$ such that
\begin{equation}\label{eq:eq60}
\mathbb E\{V(\xi_1(k+1),r_{k+1}|\xi_1(k),r_k=i\}-V(\xi_1(k),r_k=i)<-\delta\|\xi_1(k)\|^2.
\end{equation}
On the other hand, from ({\ref{eq:eq58}), we obtain
\begin{equation}\label{eq:eq61}
\lambda_{min}(\hat P_{11}(r_k))\|\xi_1(k)\|^2\leq V(\xi_1(k),r_k)\leq \lambda_{max}(\hat P_{11}(r_k))\|\xi_1(k)\|^2.
\end{equation}
Use ({\ref{eq:eq60}) and ({\ref{eq:eq61}), we have
\begin{equation}\label{eq:eq62}
\mathbb E\{V(\xi_1(k+1),r_{k+1}|\xi_1(k),r_k\}<\beta V(\xi_1(k),r_k),
\end{equation}
where
$$0<\beta=1-min_{r_k\in S}(\frac{\delta}{\lambda_{max}(\hat P_{11}(r_k))})<1.$$
By the iterative relationship, ({\ref{eq:eq62}) yields
\begin{equation}\label{eq:eq63}
\mathbb E\{V(\xi_1(k),r_k|\xi_1(0),r_0\}<\beta^k V(\xi_1(0),r_0).
\end{equation}
Therefore,
\begin{equation}\label{eq:eq64}
\mathbb E\{\|\xi_1(k)\|^2|\xi_1(0),r_0\}<\alpha\beta^k \|\xi_1(0)\|^2,
\end{equation}
where $\alpha=\frac{max_{r_k\in S}\lambda_{max}(\hat
P_{11}(r_0))}{min_{r_k\in S}\{\lambda_{min}\hat P_{11}(r_k)\}}$.
Taking  the limit on (\ref{eq:eq64}), we have
\begin{equation}\label{eq:eq65}
\lim_{k\rightarrow \infty}\mathbb E\{\|\xi_1(k)\|^2|\xi_1(0),r_0\}=0.
\end{equation}
According to Definition 6, (\ref{eq:eq55}) is asymptotically stable
in  mean square sense, and so does   system ({\ref{eq:eq42}). This
completes  the proof.

The following Theorem provides  a strict LMI condition for system
(\ref{eq:eq42}) to be mean square admissible.

 {\bf Theorem 6:} System ({\ref{eq:eq42}) is mean square admissible
if there exist positive definite matrices $P(i)>0$, $i\in S$ and a
symmetric nonsingular matrix $Q$, such that the following LMI
\begin{equation}\label{eq:eq66}
A^T(i)(\sum_{j=1}^N\lambda_{ij}P(j)+FQF^T)A(i)+C^T(i)(\sum_{j=1}^N\lambda_{ij}P(j)+FQF^T)C(i)-E^TP(i)E<0,
\end{equation}
holds,  where $F$ is a
 matrix with full column rank and satisfies $E^TF=0$.

{\bf Proof}: Let $X(i)=P(i)+FQF^T$ in (\ref{eq:eq66}), by Theorem 5,
Theorem 6 can be verified.

\section { Practical examples}
In this section, two illustrative examples related to practical
problems are proposed to demonstrate the effectiveness of our
presented approaches.

{\bf Example 1:}  Consider the modeling of oil catalytic cracking in
practical engineering (\cite{s2}, \cite{s24} ),
 which is an extremely complicated process when administration is included and its  simplified form is given as
 follows:
\begin{eqnarray}\label{eq:eq67}
\dot x_1(t)=R_{11}x_1(t)+R_{12}x_2(t)+B_1u(t)+C_1f,\nonumber\\
0=R_{21}x_1(t)+R_{22}x_2(t)+B_2u(t)+C_2f,
\end{eqnarray}
where $x_1(t)$ is a vector to be regulated   such as regenerate
temperature, valve position or  blower capacity.
 $x_2(t)$ is the vector reflecting business benefits, administration   or  policy, etc..   $u(t)$ is the regulation
 value
 and $f$ represents extra disturbances. For convenience, we consider the case of $u(t)=0, \ f=0$, then (\ref{eq:eq67}) can be
 expressed as
\begin{eqnarray}\label{eq:eq68}
E\dot{x}(t)=Rx(t),
\end{eqnarray}
where $x(t)=[x_1(t)^T\quad x_2(t)^T]^T$ is a state vector,
\begin{eqnarray}
E=\left[\begin{array}{cc}1&0\\0&0\end{array}\right], \
R=\left[\begin{array}{cc}R_{11}&R_{12}\\R_{21}&R_{22}\end{array}\right].\nonumber
\end{eqnarray}
It is obvious that (\ref{eq:eq68}) is a deterministic singular
system. However, it might happen that $R$ is subject to some random
environmental effects  (\cite{s30}-\cite{s18} ) such as
$R=A+C``noise"$. In this case, (\ref{eq:eq68}) becomes

\begin{equation}\label{eq:eq69}
\frac {Edx(t)}{dt}=Ax(t)+Cx(t)``noise".
\end{equation}
It turns out that a reasonable mathematical interpretation for the
``noise" term is the so-called white noise $\dot{w}(t)$. By
(\ref{eq:eq69}), we have
\begin{eqnarray}\label{eq:eq70}
Edx(t)=Ax(t)dt+Cx(t)dw(t).
\end{eqnarray}
In  (\ref{eq:eq70}), $A$ is called the drift matrix reflecting the
effect on the system state, while $C$ is called the diffusion matrix
reflecting the noise intensity. On the other hand, it has been
recognized that, in  many practical situations, the coefficient
matrices are not constant but random processes which can be modeled
by a Markov chain (\cite{s32}, \cite{s18} ). Therefore, we obtain
the following new equation
\begin{eqnarray}\label{eq:eq71}
Edx(t)=A(r_t)x(t)dt+C(r_t)x(t)dw(t).
\end{eqnarray}
For (\ref{eq:eq71}), if we take the following data:
\begin{eqnarray}\label{eq:eq72}
 \hspace{-5mm}&&E=\left[\begin{array}{cc}1&0\\0&0\end{array}\right],
 A(1)=\left[\begin{array}{cc}-0.5&0.7\\0.4&0.5\end{array}\right],
  A(2)=\left[\begin{array}{cc}-0.2&0.1\\0.3&0.2\end{array}\right],
C(1)=\left[\begin{array}{cc}0.4&0.2\\0&0\end{array}\right],\nonumber\\
 \hspace{-5mm}&& C(2)=\left[\begin{array}{cc}0.3&0.2\\0&0\end{array}\right],
  F=\left[\begin{array}{cc}0&1\end{array}\right]^T,
  \Pi=\left[\begin{array}{cc}-0.6&0.6\\0.5&-0.5\end{array}\right],
\end{eqnarray}
then   the solutions of LMIs (\ref{eq:eq39}) are given as follows:
\begin{eqnarray}\label{eq:eq73}
&&P(1)=\left[\begin{array}{cc}1.7492
&-0.0000\\-0.0000&1.9498\end{array}\right],
 P(2)=\left[\begin{array}{cc} 2.4364 &0.0000\\0.0000&1.9498\end{array}\right],\nonumber\\
 &&Q(1)=\left[\begin{array}{cc}-1.1149&-2.0192\end{array}\right],
  Q(2)=\left[\begin{array}{cc}-0.7417&-2.8348\end{array}\right].
\end{eqnarray}
By Theorem 3, system (3)  is mean square admissible.

{\bf Example 2:} Consider a dynamic Leontief model  of a multisector
economy without final demands (\cite{s36}, \cite{s37}, \cite{s33} )
\begin{equation}\label{eq:eq74}
x(k)=Gx(k)+E[x(k+1)-x(k)]
\end{equation}
where $x(k)$ is the vector of output levels, $G$  is the Leontief
input-output matrix, and $E$ is the capital coefficient matrix. In
economics, capital coefficient matrix $E$ is usually singular
(\cite{s36}, \cite{s37}, \cite{s33}, \cite{s34} ). However, the
above practical system are often perturbed by some environmental
noise(\cite{s17}, \cite{s18} ). Suppose that the parameter $G$ is
not completely known which is stochastically perturbed with
$G\rightarrow G+Cw(k)$, where $w(k)\in R$ is a wide stationary,
second-order process and $C$ represents the intensity of the noise.
Then this environmentally perturbed system may be described as
\begin{equation}\label{eq:eq75}
Ex(k+1)=(I-G+E)x(k)+Cx(k)w(k).
\end{equation}
Let us  further  consider another type of random fluctuation.
Suppose that the Leontief input-output matrix $G$ is  a Markov jump
process which can be modeled  by a Markov chain
(\cite{s32},\cite{s31}, \cite{s34} ), and the capital coefficient
matrix is invariant. As a result, the Leontief model is generalized
to a new model
\begin{equation}\label{eq:eq76}
Ex(k+1)=(I-G(r_t)+E)x(k)+C(r_t)x(k)w(k).
\end{equation}
Let $A(r_t)=I-G(r_t)+E$, then (\ref{eq:eq76}) is exactly the same as
(\ref{eq:eq42}). If for some concrete economical problem, the
following data are taken in (\ref{eq:eq76}):
\begin{eqnarray}
 \hspace{-5mm}&&E=\left[\begin{array}{cc}0.2&0.3\\0&0\end{array}\right],
 G(1)=\left[\begin{array}{cc}0.1&0.2\\0.3&0.1\end{array}\right],
  G(2)=\left[\begin{array}{cc}0.2&0.2\\0.4&0.5\end{array}\right],
C(1)=\left[\begin{array}{cc}0.4&-0.2\\0&0\end{array}\right],\nonumber\\
 \hspace{-5mm}&& C(2)=\left[\begin{array}{cc}0.3&-0.1\\0&0\end{array}\right],
  F=\left[\begin{array}{cc}-0.3&0.2\end{array}\right]^T,
  \Lambda=\left[\begin{array}{cc}0.4&0.6\\0.3&0.7\end{array}\right],
\end{eqnarray}
then  the solutions of LMIs (\ref{eq:eq66}) are given as follows:
\begin{eqnarray}
P(1)=\left[\begin{array}{cc}0.9548&-0.4239\\-0.4239&
0.3913\end{array}\right], P(2)=\left[\begin{array}{cc}
 0.6976&-0.3940\\-0.3940&0.3284\end{array}\right],
 Q= -11.7901.
\end{eqnarray}
Therefore, by Theorem 6, system  (46) is mean square admissible.
\section{Conclusion}
This paper has investigated the mean square admissibility of
continuous- and discrete-time singular stochastic Markov jump
systems with state-dependent noise.   Conditions for the existence
and uniqueness of the solution to the  systems considered have been
provided. In particular, new sufficient conditions for continuous-
and discrete-time singular hybrid systems with state-dependent noise
to be mean square admissible have been proposed in terms of strict
LMIs.


\begin{thebibliography}{99}



\bibitem{s35} Rosdnbrock HH. Structural properties of linear dynamical
systems. \textit{International Journal of Control} 1974;
\textbf{20}(2): 191-202.

\bibitem{s36}Luenberger DG, Arbel A.
Singular dynamic Leontief systems.\textit{ Econometrica: Journal of
the Econometric Society}  1977; 991-995.

\bibitem{s37} Luenberger DG. Dynamic equations in descriptor form.
\textit{IEEE Transactions on Automatic Control}, 1977;
\textbf{22}(3): 312-321.

\bibitem{s38} Ardema MD. \textit{Singular perturbations in systems and control}. USA: NASA, TM-62, 1983.

\bibitem{s1} Lewis FL. A survey of linear singular systems. {\em Circuits, Systems and Signal Processing}
1986; \textbf{5}(1): 3-36.

\bibitem{s2}Dai L. \textit{ Singular Control Systems}. Lecture Notes in Control
and Information Sciences. Springer: New York, 1989.

\bibitem{s39} Ailon A. On the design of output feedback for finite and infinite
pole assignment in singular systems with application to the control
problem of constrained robots. \textit{Circuits, Systems and Signal
Processing} 1994; \textbf{13}(5): 525-544.

\bibitem{s3} Xu S, Yang C. Stabilization of discrete-time singular systems: a matrix inequalities approach.
 {\em Automatica} 1999; \textbf{35}(9):1613-1617.

\bibitem{s4} Shi P,  Boukas EK, Agarwal RK. Robust control of singular continuous-time systems
 with  delays and uncertainties. \textit{ Decision and Control, Proceedings of the 39th IEEE Conference on}
 2000; \textbf{2}:1515-1520.

\bibitem{s5} Xu S, Van Dooren P,  Stefen R, Lam J. Robust stability and stabilization for singular systems
with state delay and parameter uncertainty. \textit{IEEE Transactions on Automatic Control}
2002;  \textbf{47}(7):1122-1128.

\bibitem{s6} Xu S, Lam J. Robust stability and stabilization of discrete-time singular systems:
an equivalent characterization. \textit{IEEE Transactions on Automatic Control} 2004; \textbf {49}(4): 568-574.

\bibitem{s7} Xu S, Lam J. \textit{Robust Control and Filtering of Singular
Systems}. Springer: Berlin, 2006.

\bibitem{s8}Fang C, Lee L, Chang F. Robust control analysis and design for discrete-time singular systems.
\textit{Automatica} 1994; \textbf{30}(11): 1741-1750.

\bibitem{s9} Masubuchi I, Kamitane Y, Ohara A, Suda N. $H_\infty$ control for descriptor systems:
a matrix inequalities approach. \textit{ Automatica} 1997; \textbf{3}(4): 669-673.


\bibitem{s10} Xu S, Zou Y. $H_\infty$ filtering for singular systems. \textit{IEEE Transactions on Automatic Control}
2003; \textbf{48}(12): 2217-2222.


\bibitem{s11}  Zhang L,  Huang B,  Lam J. LMI synthesis of $H_2$ and mixed $H_2/H_{\infty}$ controllers for singular
systems.
\textit{ Circuits and Systems II: Analog and Digital Signal Processing, IEEE Transactions on}
2003; \textbf{50}(9): 615-626.

\bibitem{s12} Boukas EK, Xu S, Lam J.  On stability and stabilizability of singular stochastic systems
with delays. \textit{Journal of Optimization
Theory and Applications} 2005; \textbf{127}(2): 249-262.

\bibitem{s32} Boukas EK. \textit{Control of singular systems with random abrupt
changes}. Springer Science \& Business Media, 2008.


\bibitem{s13}  Xia Y,  Zhang J, Boukas EK. Control for discrete singular hybrid systems.
\textit{ Automatica} 2008; \textbf {44}(10): 2635-2641.

\bibitem{s14} Boukas EK. On stability and stabilisation of continuous-time singular Markovian switching systems.
 \textit{ IET Control Theory and Applications} 2008; \textbf{2}(10): 884-894.

\bibitem{s15} Xia Y,  Boukas EK, Shi P, Zhang J. Stability and stabilization of continuous-time
singular hybrid systems.\textit{ Automatica} 2009; \textbf{45}(6): 1504-1509.


\bibitem{s16}  Ma S, Boukas EK, Chinniah Y. Stability and stabilization
of discrete-time singular Markov jump systems with time-varying
delay. \textit{International Journal of Robust and Nonlinear
Control}  2010; \textbf{20}(5): 531-543.

\bibitem{s30} Skorohod AV. \textit{Asymptotic Methods in the Theory of Stochastic Differential Equations}.
 American Mathematical Society Providence: RI, 1989.


\bibitem{s17} Oksendal B.\textit{Stochastic differential equations: an introduction with application}.
 Springer: New York, 1998.

\bibitem{s18} Mao X, Yuan C. \textit{Stochastic differential equations with Markovian switching}.
London: Imperial College Press, 2006.

\bibitem{s19} Zhang W,  Zhang H,  Chen B S. Generalized Lyapunov equation approach to state-dependent stochastic
stabilization/detectability criterion.
\textit{IEEE Transactions on Automatic Control}  2008;  \textbf{53}(7): 1630-1642.


\bibitem{s20} Zhang  W,  Chen B S. $\mathcal{H}$-representation and applications to generalized Lyapunov equations
 and linear stochastic systems. \textit{IEEE Transactions on Automatic Control}  2012;  \textbf{57}(12): 3009-3022.

\bibitem{s21}Boukas  EK. Stabilization of stochastic singular nonlinear hybrid systems.
 \textit{Nonlinear Analysis}  2006; \textbf{64}(2): 217-228.

\bibitem{s22} Zhang Q, Xing S. Stability analysis and optimal control of stochastic singular systems.
\textit{Optimization Letters} 2014; \textbf{8}(6): 1905-1920.

\bibitem{s23} Han C, Wu L, Shi P, Zheng Q. Passivity and passification
of T-S fuzzy descriptor systems with stochastic perturbation and
time delay. \textit{IET Control Theory and Applications} 2013; \textbf{7}(13): 1711-1724.


\bibitem{s24} Zhang W, Zhao, Y, Sheng L. Some remarks on stability of stochastic singular
systems with state-dependent noise. \textit{ Automatica}  2015;
\textbf{51}: 273-277.


\bibitem{s25} Ho D, Shi X, Wang Z, Gao Z. Filtering for a
class of stochastic descriptor systems. \textit{In Proceedings of
the International Conference on Dynamics of Continuous, Discrete and
Impulsive Systems, Canada}  2005; \textbf{2}: 848-853.

\bibitem{s26} Gao  Z, Shi X. Observer-based controller design for stochastic descriptor systems with
Brownian motions. \textit{Automatica} 2013; \textbf{49}(7): 2229-2235.

\bibitem{s27} Huang  L,  Mao X. Stability of singular stochastic systems with Markovian
switching. \textit{IEEE Transactions on Automatic Control}  2011; \textbf{56}(2): 424-429.

\bibitem{s28} Xia  J. \textit{Robust Control and Filter for Continuous Stochastic Time-Delay  Systems}.
Ph.D. Dissertation, Nanjing University of Science and Technology,
China, 2007.

\bibitem{s29} Bellman RE. \textit{Introduction to Matrix Analysis}. SIAM, 1995.


\bibitem{s31} West GR.  A stochastic analysis of an input-output model. \textit{Econometrica: Journal of the
Econometric Society}  1986; \textbf{54}(2): 363-374.



\bibitem{s33} Mao W J. An LMI
approach to D-stability and D-stabilization of linear discrete
singular systems with state delay. \textit{Applied Mathematics and
Computation} 2011; \textbf{218}(5): 1694-1704.

\bibitem{s34} Wu X, Jiang L. Computer analysis algorithm for stability of the
extended dynamic Leontief input-output model. \textit{ Proceedings
of International Conference on Computational Intelligent and Natural
Computing, Wuhan, China} 2009; \textbf{2}: 379-382.

\end{thebibliography}
\end{document}